\documentclass{amsart}

\usepackage{fancyhdr}
\usepackage{lastpage}
\usepackage{stmaryrd,yhmath}

\newcommand{\bdis}{\begin{displaymath}}
\newcommand{\edis}{\end{displaymath}}
\newcommand{\be}{\begin{equation}}
\newcommand{\ee}{\end{equation}}
\newcommand{\mbb}{\mathbb}
\newcommand{\mcal}{\mathcal}

\newcommand{\vp}{\varphi}

\newcommand{\vth}{\vartheta}

\newcommand{\zf}{\zeta\left(\frac{1}{2}+it\right)}

\newcommand{\okT}{\overset{k}{T}}
\newcommand{\onT}{\overset{0}{T}}

\theoremstyle{definition}

\theoremstyle{remark}
\newtheorem{remark}[]{Remark}

\newtheorem*{mydef11}{{\bf Theorem 1}}

\newtheorem*{mydef12}{{\bf Theorem 2}}

\newtheorem*{mydef13}{{\bf Theorem 3}}

\newtheorem*{mydef41}{{\bf Corollary 1}}

\newtheorem*{mydef42}{{\bf Corollary 2}}

\newtheorem*{mydef6}{{\bf Example}}

\numberwithin{equation}{section}

\begin{document}

\title{Jacob's ladders and laws that control chaotic behavior of the measures of reversely iterated segments}

\author{Jan Moser}

\address{Department of Mathematical Analysis and Numerical Mathematics, Comenius University, Mlynska Dolina M105, 842 48 Bratislava, SLOVAKIA}

\email{jan.mozer@fmph.uniba.sk}


\begin{abstract}
The main subject to study in this paper are properties of the sequence of reversely iterated segments. Especially, we will examine
properties of chaotic behavior of the sequence of measures of corresponding segments. Our results are not accessible within current
methods in the theory of Riemann zeta-function.
\end{abstract}
\maketitle

\section{Introduction}

\subsection{}

Let us start with some notions and formulae to be reminded:
\begin{itemize}
\item[(A)] the sequence
\bdis
\{\okT\}_{k=1}^{k_0}
\edis
is defined by (see \cite{3}, (5.1))
\bdis
\vp_1(\okT)=\overset{k-1}{T},\ k=1,\dots,k_0,\ \onT=T,\ T\geq T_0[\vp_1] ,
\edis
where $k_0\in\mbb{N}$ is an arbitrary fixed number and $\vp_1(t)$ is the Jacob's ladder;

\item[(B)] next (see \cite{3}, (1.3))
\be \label{1.1}
\begin{split}
& \tilde{Z}^2(t)=\frac{{\rm d}\vp_1(t)}{{\rm d}t}=\frac{Z^2(t)}{2\Phi'_\vp[\vp(t)]}=\frac{\left|\zf\right|^2}{\omega(t)}, \\
& \omega(t)=\left\{ 1+\mcal{O}\left(\frac{\ln\ln t}{\ln t}\right)\right\}\ln t.
\end{split}
\ee
where
\be \label{1.2}
\begin{split}
 & Z(t)=e^{i\vth(t)}\zf, \\
 & \vth(t)=-\frac t2\ln \pi +\text{Im}\ln\Gamma\left(\frac 14+i\frac t2\right),
\end{split}
\ee
\end{itemize}

\subsection{}

We have proved the following theorem (see \cite{3}, (2.1) -- (2.7)): for every $L_2$-orthogonal system
\bdis
\{ f_n(t)\}_{n=1}^\infty,\ t\in [0,2l],\ l=o\left(\frac{T}{\ln T}\right),\ T\to\infty
\edis
there is a continuum set of $L_2$-orthogonal systems
\bdis
\begin{split}
& \{F_n(t;T,k,l)\}_{n=1}^\infty= \\
& = \left\{ f_n(\vp_1(t)-T)\prod_{r=0}^{k-1}\left|\tilde{Z}[\vp_1^r(t)]\right|\right\}_{n=1}^\infty, \
t\in[\okT,\overset{k}{\wideparen{T+2l}}],
\end{split}
\edis
where
\be \label{1.3}
\begin{split}
& \vp_1\{[\okT,\overset{k}{\wideparen{T+2l}}]\}=[\overset{k-1}{T},\overset{k-1}{\wideparen{T+2l}}],\ k=1,\dots,k_0, \\
& [\onT,\overset{0}{\wideparen{T+2l}}]=[T,T+2l],
\end{split}
\ee
i.e. the following formula is valid
\bdis
\begin{split}
& \int_{\overset{k}{T}}^{\overset{k}{\wideparen{T+2l}}}f_m(\vp_1^k(t)-T)f_n(\vp_1^k(t)-T)\prod_{r=0}^{k-1}
\tilde{Z}^2[\vp_1^r(t)]{\rm d}t=\\
& = \left\{\begin{array}{rcl} 0 & , & m\not= n, \\ A_n & , & m=n,  \end{array}\right.\quad
A_n=\int_0^{2l} f_n^2(t){\rm d}t.
\end{split}
\edis

\begin{remark}
It is clear that the base of above mentioned result is new notion of reverse iterations (comp. (1.3)) in the theory of Riemann
$\zf$-function.
\end{remark}

In this paper we will study the sequence of reverse iterations
\bdis
\left\{ [\overset{r}{T},\overset{r}{\wideparen{T+H}}]\right\}_{r=0}^k,\ k=1,\dots,k_0
\edis
alone. Namely, we will focus on properties of the sequence of real numbers (measures of corresponding segments)
\bdis
\left\{ |[\overset{r}{T},\overset{r}{\wideparen{T+H}}]|\right\}_{r=0}^k.
\edis

\begin{remark}
Results of this paper are not accessible by current methods of the theory of Riemann zeta-function. We mention explicitly that our results
are valid also in the microscopic case
\bdis
H\in \left(\left. 0,\frac{A}{\ln T}\right]\right., \quad T\to\infty.
\edis
\end{remark}

\section{Theorem 1 and motivation behind it}

\subsection{}

Let us remind that the segments
\bdis
[\overset{r}{T},\overset{r}{\wideparen{T+H}}],\ r=0,1,\dots,k
\edis
are components of disconnected set
\be \label{2.1}
\Delta(T,H,k)=\bigcup_{r=0}^k [\overset{r}{T},\overset{r}{\wideparen{T+H}}],\ k=1,\dots,k_0,
\ee
(comp. \cite{3}, (2.9)). Properties of the set (2.1) are listed below (see \cite{3}, (2.5) -- (2.7)):
\be \label{2.2}
\begin{split}
& H=o\left(\frac{T}{\ln T}\right) \ \Rightarrow \\
& |[\overset{k}{T},\overset{k}{\wideparen{T+H}}]|=\overset{k}{\wideparen{T+H}}-\okT=o\left(\frac{T}{\ln T}\right),
\end{split}
\ee
\be \label{2.3}
|[\overset{k-1}{\wideparen{T+H}},\okT]|=\okT-\overset{k-1}{\wideparen{T+H}}\sim (1-c)\pi(T);\ \pi(T)\sim\frac{T}{\ln T},
\ee
\be \label{2.4}
[T,T+H]\prec [\overset{1}{T},\overset{1}{\wideparen{T+H}}]\prec \dots \prec [\overset{k}{T},\overset{k}{\wideparen{T+H}}]\prec \dots ,
\ee
where $c$ is the Euler's constant and $\pi(T)$ is the prime-counting function.

\begin{remark}
Consequently, the asymptotic behavior of our disconnected set (2.1) is as follows (see (2.2), (2.3)): if $T\to\infty$ then the components of the set
(2.1) recede unboundedly each from other and all together are receding to infinity. Hence, the set (2.1) behaves as a kind of one-dimensional
Friedmann-Hubble expanding universe.
\end{remark}

Furthermore, we notice explicitly that the distance $\rho_l$ of the two consecutive segments
\bdis
[\overset{l-1}{T},\overset{l-1}{\wideparen{T+H}}],\ [\overset{l}{T},\overset{l}{\wideparen{T+H}}],\quad l=1,2,\dots,k
\edis
is extremely big one, namely (see (2.3))
\be \label{2.5}
\rho_l\sim (1-c)\frac{T}{\ln T}\to\infty,\quad T\to\infty.
\ee

\begin{remark}
Since the sequence
\bdis
\{[\overset{r}{T},\overset{r}{\wideparen{T+H}}]\}_{r=0}^k
\edis
is extremely sparse one (see (2.5)) then we may assume that the behavior of the measures
\bdis
\{|[\overset{r}{T},\overset{r}{\wideparen{T+H}}]|\}_{r=0}^k
\edis
is chaotic one.
\end{remark}

Consequently, in correspondence with Remark 3, we wish to obtain some law controlling this chaotic
behavior. In this direction, the following theorem holds true.

\begin{mydef11}
Let
\be \label{2.6}
1\leq n\leq k_0,\quad \bar{H}=o\left(\frac{T}{\ln T}\right),
\ee
and let the inequality
\be \label{2.7}
|[\overset{n}{T},\overset{n}{\wideparen{T+\bar{H}}}]|=
\overset{n}{\wideparen{T+\bar{H}}}-\overset{n}{T}\geq T^{1/3+\epsilon},\ T\to\infty
\ee
hold true for $\epsilon>0$ - an arbitrary small fixed number. Then we have that
\be \label{2.8}
\overset{n}{\wideparen{T+\bar{H}}}-\overset{n}{T}\sim \bar{H},\ T\to\infty.
\ee
\end{mydef11}

\subsection{}

Next, let us remind
\begin{itemize}
 \item[(A)] the Hardy-Littlewood-Ingham formula
 \be \label{2.9}
 \begin{split}
  & \int_0^T\left|\zf\right|^2{\rm d}t=
  T\ln T+(2c-1-\ln 2\pi)T+R(T),
 \end{split}
 \ee
 with the Balasubramanian's estimate (for example)
 \be \label{2.10}
 R(T)=\mcal{O}(T^{1/3})
 \ee
 of the error term in (2.9);

 \item[(B)] the Good's $\Omega$-theorem that states
 \be \label{2.11}
 R(T)=\Omega(T^{1/4}),\ T\to\infty ;
 \ee

 \item[(C)] our almost exact formula (see \cite{1}, (2.1), (2.2),
 $\frac y2\to \vp_1(t)$)
 \be \label{2.12}
 \begin{split}
  & \int_0^T\left|\zf\right|^2{\rm d}t=\\
  & = \vp_1(T)\ln\vp_1(T)+(c-\ln 2\pi)\vp_1(T)+c_0+\mcal{O}\left(\frac{\ln T}{T}\right),\
  T\to\infty,
 \end{split}
\ee
where $c$ is the Euler's constant and $c_0$ is the constant from the
Titchmarsh-Kober-Atkinson formula (see \cite{4}, p. 141).
\end{itemize}

Our discussion concerning formulae (2.10) -- (2.12) see in \cite{1}, pp.
416, 417.

\begin{remark}
 Consequently, we have proved in \cite{1} that classical Hardy-Littlewood integral (1918)
\bdis
\int_0^T\left|\zf\right|^2{\rm d}t
\edis
has -- in addition to the Hardy-Littlewood (and other similar) expressions possessing
unbounded errors (as $T\to\infty$), (comp. (2.10), (2.11)) -- infinite set of almost
exact expressions (2.12).
\end{remark}

\begin{remark}
It is clear -- in context of (2.10), (2.11) -- that our Theorem 1 will be true for every
improvement of the exponent $\frac 13$:
\bdis
\frac 13 \longrightarrow a\in \left(\frac 14,\frac 13\right).
\edis
\end{remark}

\section{Proof of Theorem 1}

First of all, it follows from (2.9) and (2.10) that
\be \label{3.1}
\begin{split}
 & \int_T^{T+U}\left|\zf\right|^2{\rm d}t \sim U\ln T,\ T\to\infty, \\
 & T^{1/3+\epsilon}\leq U=o\left(\frac{T}{\ln T}\right).
\end{split}
\ee
Next, we use, together with (3.1), our formula
\be \label{3.2}
\int_{\overset{k}{T}}^{\overset{k}{\wideparen{T+H}}}\left|\zf\right|^2{\rm d}t\sim
(\overset{k-1}{\wideparen{T+H}}-\overset{k-1}{T})\ln T
\ee
that follows from \cite{3}, (1.1) -- (1.3), (7.4) with
\bdis
[T,T+H]\longrightarrow [\overset{k-1}{T},\overset{k-1}{\wideparen{T+H}}],\
[\overset{1}{T},\overset{1}{\wideparen{T+H}}]\longrightarrow
[\overset{k}{T},\overset{k}{\wideparen{T+H}}].
\edis
Of course, we have (see (2.2))
\be \label{3.3}
\begin{split}
 & H=o\left(\frac{T}{\ln T}\right) \ \Rightarrow \\
 & \Rightarrow\
 \overset{k}{\wideparen{T+H}}-\overset{k}{T}=o\left(\frac{T}{\ln T}\right),\quad
 T\to\infty,\ k=1,\dots,k_0.
\end{split}
\ee
Further, if $n,\bar{H}$ fulfill the conditions (2.6) and (2.7) then we have
(see (3.1), (3.2)) that
\bdis
\int_{\overset{n}{T}}^{\overset{n}{\wideparen{T+\bar{H}}}}\left|\zf\right|^2{\rm d}t\sim
(\overset{n}{\wideparen{T+\bar{H}}}-\overset{n}{T})\ln T,
\edis
and
\bdis
\int_{\overset{n}{T}}^{\overset{n}{\wideparen{T+\bar{H}}}}\left|\zf\right|^2{\rm d}t\sim
(\overset{n-1}{\wideparen{T+\bar{H}}}-\overset{n-1}{T})\ln T,
\edis
i.e.
\bdis
\overset{n}{\wideparen{T+\bar{H}}}-\overset{n}{T}\sim
\overset{n-1}{\wideparen{T+\bar{H}}}-\overset{n-1}{T},\quad T\to\infty,
\edis
and, consequently,
\bdis
\overset{n-1}{\wideparen{T+\bar{H}}}-\overset{n-1}{T}\sim
\overset{n-2}{\wideparen{T+\bar{H}}}-\overset{n-2}{T}\sim \dots \sim
T+\bar{H}-T=\bar{H},
\edis
(see also (3.3)). Thus, we have that
\bdis
\overset{n}{\wideparen{T+\bar{H}}}-\overset{n}{T}\sim \bar{H},\quad T\to\infty,
\edis
i.e. the assertion (2.8) is verified.

\section{Consequences of Theorem 1}

\subsection{}

\begin{mydef41}
Let
\be \label{4.1}
H_1=A(T)T^{1/3+\epsilon},\quad 0<A(T)<1,
\ee
for example
\bdis
H_1=\frac 12T^{1/3+\epsilon},\ \frac{1}{\ln\ln T}T^{1/3+e\epsilon},\dots
\edis
Then
\be \label{4.2}
\overset{k}{\wideparen{T+H_1}}-\overset{k}{T}< T^{1/3+e\epsilon},\quad
k=1,\dots,k_0.
\ee
\end{mydef41}

\begin{remark}
Hence, in the case (4.1) we have that all members of the sequence
\bdis
\left\{|[\overset{k}{T},\overset{k}{\wideparen{T+H_1}}]|\right\}_{k=1}^{k_0}
\edis
are lying below the level $T^{1/3+\epsilon}$, (see (4.2)).
\end{remark}

\subsection{}

Next, as a consequence of Corollary 1, we have
\begin{mydef42}
If
\be \label{4.5}
H_2=B(T)T^{1/3+\epsilon},\quad B(T)>1,
\ee
for example
\bdis
H_2= 2T^{1/3+\epsilon},\ T^{1/3+\epsilon}\ln T,\dots
\edis
and there is some
\bdis
n:\ 1\leq n<k_0
\edis
such that
\be \label{4.4}
\overset{n}{\wideparen{T+H_2}}-\overset{n}{T}<A(T)T^{1/3+\epsilon}
\ee
(see (4.1)), then
\be \label{4.5}
\overset{k}{\wideparen{T+H_2}}-\okT< T^{1/3+\epsilon},\quad k=n+1,\dots,k_0.
\ee
\end{mydef42}

\begin{remark}
Consequently, in the case (4.3), (4.4) the second jump of the sequence
\bdis
\left\{|[\okT,\overset{k}{\wideparen{T+H_2}}]|\right\}_{k=1}^{k_0}
\edis
over the segment
\be \label{4.6}
[(1-\epsilon)T^{1/3+\epsilon},(1+\epsilon)T^{1/3+\epsilon}]
\ee
is forbidden. In other words, the oscillations of the sequence of measures about the measure of the segment (4.6) are forbidden.
\end{remark}

\section{An estimate from below}

\subsection{}

We will use the following in this section:
\begin{itemize}
\item[(A)] the estimate
\bdis
\begin{split}
& H=o\left(\frac{T}{\ln T}\right) \ \Rightarrow \\
& \int_{\overset{k}{T}}^{\overset{k}{\wideparen{T+H}}}\left|\zf\right|^2{\rm d}t> (1-\epsilon)(\overset{k-1}{\wideparen{T+H}}-\overset{k}{T})\ln T,\
k=1,\dots,k_0,
\end{split}
\edis
that follows from the asymptotic formula (3.2), i.e. we have that
\be \label{5.1}
\begin{split}
& \int_{\overset{1}{T}}^{\overset{1}{\wideparen{T+H}}}\left|\zf\right|^2{\rm d}t>(1-\epsilon)H\ln T, \\
& \int_{\overset{2}{T}}^{\overset{2}{\wideparen{T+H}}}\left|\zf\right|^2{\rm d}t> (1-\epsilon)(\overset{1}{\wideparen{T+H}}-\overset{1}{T})\ln T, \\
& \vdots \\
&\int_{\overset{k_0}{T}}^{\overset{k_0}{\wideparen{T+H}}}\left|\zf\right|^2{\rm d}t> (1-\epsilon)(\overset{k_0-1}{\wideparen{T+H}}-\overset{k_0}{T})\ln T;
\end{split}
\ee

\item[(B)] the property (see (2.1) and \cite{3}, sec. 4.1)
\be \label{5.2}
\tau\in\Delta(T,H,k) \ \Rightarrow \
\tau\in \left[ T,T+\mcal{O}\left(\frac{T}{\ln T}\right)\right]\subset [T,2T].
\ee
\end{itemize}

\subsection{}

Since (comp. \cite{4}, p. 99)
\be \label{5.3}
\begin{split}
& \left|\zf\right|<t^{1/6},\ t\to\infty \ \Rightarrow \\
& \left|\zf\right|^2<t^{1/3},\ t\to\infty,\ t\to\infty,
\end{split}
\ee
then we have (see (5.2), (5.3)), that
\be \label{5.4}
\begin{split}
& t\in \Delta(T,H,k) \ \Rightarrow \\
& \left|\zf\right|^2<\sqrt[3]{2}T^{1/3}<2T^{1/3},\ T\to\infty,
\end{split}
\ee
without any hypothesis. Consequently, we have (see (5.1), (5.4)) following estimates
\be \label{5.5}
\begin{split}
& |[\overset{1}{T},\overset{1}{\wideparen{T+H}}]|>\frac{1-\epsilon}{2}HT^{-1/3}\ln T, \\
& |[\overset{2}{T},\overset{2}{\wideparen{T+H}}]|>\left(\frac{1-\epsilon}{2}HT^{-1/3}\ln T\right)^2H, \\
& \vdots \\
& |[\overset{k_0}{T},\overset{k_0}{\wideparen{T+H}}]|>\left(\frac{1-\epsilon}{2}HT^{-1/3}\ln T\right)^{k_0}H>
\left(\frac 14 T^{-1/3}\ln T\right)^{k_0}H=\\
& = \left(\frac{\ln^3T}{64T}\right)^{k_0/3}H,\quad \epsilon\in (0,1/2).
\end{split}
\ee
Since
\bdis
0<\frac{\ln T}{4T^{1/3}}<1,\ T\to\infty,
\edis
then we have the following

\begin{mydef12}
\be \label{5.6}
\begin{split}
& H=o\left(\frac{T}{\ln T}\right) \ \Rightarrow \\
& |[\overset{k}{T},\overset{k}{\wideparen{T+H}}]|>\left(\frac{\ln^3T}{64T}\right)^{k_0/3}H,\ k=1,\dots,k_0,\ T\to\infty.
\end{split}
\ee
\end{mydef12}

\begin{mydef6}
If
\bdis
H=1,\ k_0=3000
\edis
then (see (5.6))
\bdis
|[\overset{k}{T},\overset{k}{\wideparen{T+H}}]|>\left(\frac{\ln T}{64T}\right)^{1000},\ k=1,\dots,3000.
\edis
\end{mydef6}

\begin{remark}
It appears that only advantage of the estimate (5.6) is, probably, its non-triviality.
\end{remark}

\section{Riemann hypothesis and our estimate from below}

The following estimate
\bdis
\left|\zf\right|< Be^{A\frac{\ln t}{\ln\ln t}},\quad t\to\infty
\edis
holds true on the Riemann hypothesis (see \cite{4}, p. 300). We use this estimate in the form
\bdis
\left|\zf\right|<t^{\frac{C}{\ln\ln t}},\quad t\to\infty.
\edis
Thus we have (comp. (5.2), (5.4)) that
\be \label{6.1}
\begin{split}
 & t\in\Delta(T,H,k) \ \Rightarrow \
 \left|\zf\right|^2< (2T)^{\frac{2C}{\ln\ln(2 T)}}<(2T)^{\frac{2C}{\ln\ln T}}< \\
 & < 2^{\frac{2C}{\ln\ln T}}T^{\frac{2C}{\ln\ln T}}<
 (1+\epsilon)T^{\frac{2C}{\ln\ln T}},\quad \epsilon\in (0,1/2),\quad T\to\infty.
\end{split}
\ee
Now, we obtain from (5.1), (6.1), (comp. (5.2), (5.5)) that
\bdis
\begin{split}
 & |[\overset{k_0}{T},\overset{k_0}{\wideparen{T+H}}]|>\left(\frac{1-\epsilon}{1+\epsilon}\right)^{k_0}
 HT^{-k_0\frac{2C}{\ln\ln T}}\ln^{k_0}T>\\
 & > \frac{1}{3^{k_0}}H T^{-k_0\frac{2C}{\ln\ln T}}\ln^{k_0}T= \\
 & = T^{-k_0\frac{\ln 3}{\ln T}}H T^{-k_0\frac{2C}{\ln\ln T}}T^{k_0\frac{\ln\ln T}{\ln T}}>
 H T^{-k_0\frac{2D}{\ln\ln T}}.
\end{split}
\edis
Hence, the following theorem holds true.

\begin{mydef13}
On Riemann hypothesis we have
\be \label{6.2}
\begin{split}
 & H=o\left(\frac{T}{\ln T}\right) \ \Rightarrow \\
 & |[\overset{k}{T},\overset{k}{\wideparen{T+H}}]|>
 H T^{-k_0\frac{2D}{\ln\ln T}},\quad k=1,\dots,k_0,\quad T\to\infty.
\end{split}
\ee
\end{mydef13}

\begin{remark}
The conditional estimate (6.2) is effective particularly in the case
\be \label{6.3}
H=T^\Delta,\ 0<\Delta<1.
\ee
Namely, in this case we obtain from (6.2)
\be \label{6.4}
|[\overset{k}{T},\overset{k}{\wideparen{T+H}}]|>T^{\Delta-o(1)},\quad T\to\infty.
\ee
\end{remark}

\begin{remark}
It was expected that the Riemann hypothesis has essential influence on that estimate from below. Actually, we have in the case (6.3) that:
\begin{itemize}
 \item[(A)] without any hypothesis (see (5.6))
 \be \label{6.5}
 |[\overset{k}{T},\overset{k}{\wideparen{T+T^\Delta}}]|>
 \left(\frac 14\ln T\right)^{k_0}T^{\Delta-\frac{k_0}{3}},
 \ee
 where
 \be \label{6.6}
 k_0\geq 3 \ \Rightarrow \ \Delta-\frac{k_0}{3}<0;
 \ee

 \item[(B)] on the Riemann hypothesis (see (6.4))
 \be \label{6.7}
 |[\overset{k}{T},\overset{k}{\wideparen{T+T^\Delta}}]|>T^{\Delta-o(1)},
 \ee
 where
 \be \label{6.8}
 0<\Delta-o(1)\to\Delta \ \text{as}\ T\to\infty.
 \ee
\end{itemize}
\end{remark}

\begin{mydef6}
In the case
\bdis
\Delta=\frac 13+\epsilon
\edis
(see Theorem 1) we have, on Riemann hypothesis, that (see (6.4))
\bdis
|[\overset{k}{T},\overset{k}{\wideparen{T+T^{1/3+\epsilon}}}]|>
T^{1/3+\epsilon-o(1)},\quad k=1,\dots,k_0,\quad T\to\infty.
\edis
\end{mydef6}

\thanks{I would like to thank Michal Demetrian for his help with electronic version of this paper.}

\end{document}